\documentclass[12pt,french]{amsart}

\usepackage{amssymb}
\usepackage{graphicx}
\usepackage{epsf}
\usepackage{amsmath}
\usepackage[latin1]{inputenc}
\usepackage{amsfonts}
\usepackage{color}

\newtheorem{theo}{Theorem}
\newtheorem{prop}{Proposition}

\newtheorem{defi}[theo]{Definition}

\newtheorem{rema}[theo]{Remark}

\newenvironment{resume}{\footnotesize\quotation}

\def\cf{{\it cf. }}
\def\ie{{\it i.e. }}

\def\N{{\mathbb N}}

\def\S{{\mathcal S}}

\def\sha{{\,\#\,}}

\def\shuffle{{\,\raise
1pt\hbox{$\scriptscriptstyle\cup{\mskip-4mu}\cup$}\,}}

\def\Young#1{\vbox{\smallskip\offinterlineskip
    \halign{&\vbox{##}\kern-\Thickness\cr #1}}}

\newdimen\Squaresize \Squaresize=12pt
\newdimen\Thickness \Thickness=.1pt
\newdimen\Correction \Correction=7pt

\def\Vide#1{\hbox{
       \vbox to \Squaresize{\vss
          \hbox to \Squaresize{\hss#1 \hss}\vss}
    \hskip-\Correction}
   \kern-\Thickness}

\def\Carre#1{\hbox{\vrule width \Thickness
   \vbox to \Squaresize{\hrule height \Thickness\vss
      \hbox to \Squaresize{\hss$\scriptstyle#1$\hss}
   \vss\hrule height\Thickness}
   \unskip\vrule width \Thickness}
   \kern-\Thickness}

\title[Product of trees and Catalan alternative tableaux]{The product of trees in the Loday-Ronco algebra through Catalan alternative tableaux}

\author{J.-C.~Aval, X. Viennot}
\address[Jean-Christophe Aval]{LaBRI\\ Universit\'e Bordeaux 1\\ 351 cours
 de la Lib\'eration\\ 33405 Talence cedex\\ FRANCE}
\email{aval@labri.fr}
\urladdr{http://www.labri.fr/perso/aval}
\address[Xavier Viennot]{LaBRI\\ Universit\'e Bordeaux 1\\ 351 cours
 de la Lib\'eration\\ 33405 Talence cedex\\ FRANCE}
\email{viennot@labri.fr}
\urladdr{http://www.labri.fr/perso/viennot}

\thanks{This research has been supported by the ANR (project MARS/06-BLAN-0193)}

\begin{document} 

\maketitle 

\begin{resume}
The aim of this note is to show how the introduction of certain tableaux, called {\em Catalan alternative tableaux}, provides a very simple and elegant description of the product in the Hopf algebra of binary trees defined by Loday and Ronco. Moreover, we use this description to introduce a new associative product on the space of binary trees.
\end{resume}

\bigskip
\bigskip

\section{Introduction}

Loday and Ronco defined in \cite{LR} an interesting Hopf algebra structure on the linear span of rooted planar binary trees. This algebra is defined as a sub-algebra of the Malvenuto-Reutenauer Hopf algebra of permutations. 
Let $S_n$ be the symmetric group and $k$ be a ground field. We denote by $k[S_n]$ the group algebra. Malvenuto and  Reutenauer construct in \cite{MR} a Hopf algebra structure on
$$k[S_\infty]=\bigoplus_{n\ge 0}k[S_n].$$
It is worth to recall here that the Malvenuto-Reutenauer algebra contains the sum of Solomon descent algebras $Sol_\infty=\bigoplus_{n\ge 0}Sol_n$ with $Sol_n$ of dimension $2^{n-1}$. 

In \cite{LR}, Loday and Ronco define a sub-Hopf algebra of $k[S_\infty]$:
$$k[Y_\infty]=\bigoplus_{n\ge 0}k[Y_n]$$
where $Y_n$ is the set of planar binary trees with $n$ internal vertices.

The aim of this work is to present a very simple presentation for the product of two trees in $k[Y_\infty]$ through the use of {\em Catalan alternative tableaux}. These objects were introduced by X. Viennot \cite{AT} as a special case of alternative tableaux, which are in bijection with permutation tableaux.
Permutations tableaux were introduced by  E. Steingrimsson and L. Williams \cite{SW}, as a subclass of $\Gamma$-diagram defined by A. Postnikov \cite{post}.  This notion was used by S. Corteel and L. Williams \cite{CW1,CW2} in the study of the physical model named PASEP (partially asymmetric exclusion process), see for example the seminal paper by B. Derrida and al. \cite{PASEP}. These tableaux are also related to the study of total positivity for Grassmannian \cite{williams}. Both permutation and alternative tableaux are in bijection with permutations, see for example P. Nadeau \cite{nadeau}. The advantage of alternative tableaux is to preserve the symmetry between rows and columns.

This new interpretation of the Loday-Ronco product motivates the introduction of a new associative product, that we call the $\sha$ product, on the space of binary trees. This new product is studied in \cite{chap} by F. Chapoton.

This article is constructed as follows: in Section 2 we recall the definition of the Loday-Ronco algebra, then we introduce in Section 3 the Catalan alternative tableaux and prove that they are in bijection with binary trees; we state and prove in Section 4 the main result of this work, and we introduce the new $\sha$ product in Section 5.

\section{The Loday-Ronco Hopf algebra}

We recall the definition of the Loday-Ronco product of binary trees. Since this product  is inherited from the Malvenuto-Reutenauer product of permutations, we shall first recall the definition of the product in $k[S_\infty]$, denoted by $*$. We refer to \cite{MR} for more details and only recall briefly the definition. 

Let $u=u_1\,u_2\,\dots,u_k$ be a $k$-tuple of distinct integers. We define the {\em standardization} of $u$ and denote it by Std$(u$) as the unique permutation $\sigma\in S_k$ that preserves the relative order of the $u_i$'s, {\ie}
$$\sigma_i<\sigma_j \Longleftrightarrow u_i<u_j.$$ 

For example, Std$(3275)=2143$. Conversely, for $\sigma\in\S_k$ a permutation and $A=\{a_1,a_2,\dots,a_k\}$ a set of $k$ (distinct) integers, we define $\sigma_{|A}$ the $k$-tuple with distinct entries in $A$ such that Std$(\sigma_{|A})=\sigma$. With this notation we may define the product $*$ in $k[S_\infty]$ as follows. Let $\sigma\in S_k$ and $\tau\in S_l$. We set
$$\sigma*\tau=\sum_{A\sqcup B=\{1,2,\dots,k+l\}}\sigma_{|A}.\tau_{|B}$$
where $\sqcup$ denotes the disjoint union, and $.$ stands for concatenation.

For example :
$$12*213=12\,435+13\,425+14\,325+15\,324+23\,415+24\,315+25\,314+34\,215+35\,214+45\,213.$$

\begin{rema}
The product $*$ that we consider is sometimes known as the product in the {\em dual} Malvenuto-Reutenauer algebra. But it is the one used in \cite{LR} to define the Loday-Ronco algebra, that we shall now describe.
\end{rema}

\medskip

Let $Y_n$ denote the set of binary trees with $n$ internal vertices. We recall that the cardinality of $Y_n$ is given by the $n$-th Catalan number $C_n=\frac 1 {n+1} {2n \choose n}$.

Let $\tilde Y_n$ denote the set of {\em increasing binary trees}, \ie of binary trees such that each internal vertex has a distinct label in $\{1,\dots,n\}$, and such that the labels increase along the tree.

It is well known that increasing binary trees are in bijection with permutations: 
to obtain the permutation from the tree, you just have to read the labels from left to right.

Below is an example of a plane binary tree with 8 internal vertices, with a increasing binary tree with the same underlying tree, and with the corresponding permutation $\sigma\in S_n$.

\centerline{\epsffile{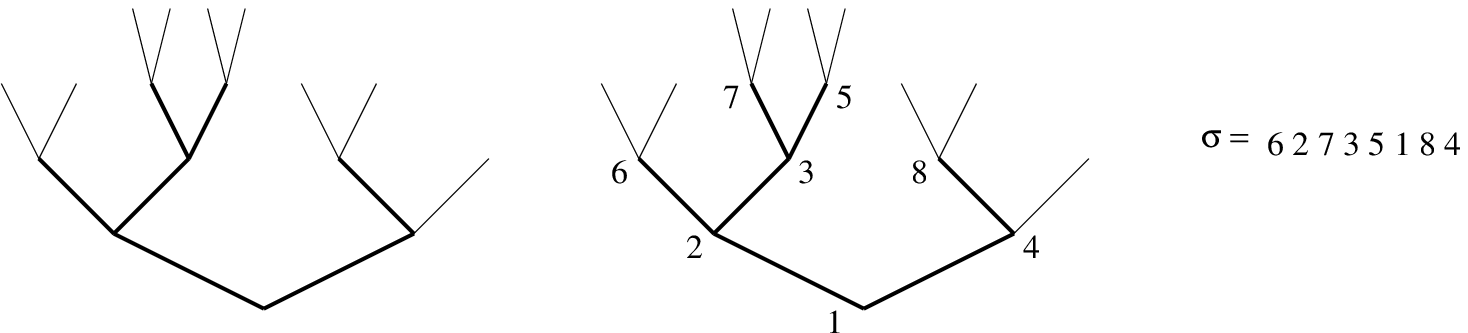}}

We denote by $\Psi: S_n \rightarrow Y_n$ the composition of the bijection $S_n\simeq \tilde Y_n$ with the projection $\tilde Y_n \rightarrow Y_n$ which consists in forgetting the labels.
The induced linear map $\Psi: k[S_n] \rightarrow k[Y_n]$ has a linear dual $\Psi^*: k[Y_n] \rightarrow k[S_n]$ obtained by identifying each basis with its own dual. For example
$$\Psi^*\Big(\epsffile{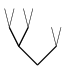}\Big)=3412+4312+2413+4213+2314+3214.$$

We also define for any tree $T$ the set 
$$Z_T=\{\sigma\in S_n\ /\ \Psi(\sigma)=T\}$$
so that $\Psi^*(T)=\sum_{\sigma\in Z_T} \sigma$.

The inclusion map $\Psi^*$ gives rise to a graded linear map $\Psi^*:k[Y_\infty]\rightarrow k[S_\infty]$ and the main result in the construction of the Loday-Ronco algebra may now be stated as (Theorem 3.1 in \cite{LR}):

\begin{theo}
The image of the inclusion map $\Psi^*:k[Y_\infty]\rightarrow k[S_\infty]$ is a sub-Hopf algebra of $k[S_\infty]$. So, $k[Y_\infty]$ inherits a structure of Hopf algebra.
\end{theo}

\section{Trees and Catalan alternative tableaux}

We now present the Catalan alternative tableaux. 
Let us denote by $\N$ the set of nonnegative integers. A Catalan alternative tableau in given by
\begin{itemize}
\item a path in $\N\times \N$ from $\{0\}\times \N$ to $\N\times \{0\}$ made of $(0,1)$ and $(1,0)$ steps. The length of the path is called the {\em size} of the tableau, and the cells below the path are simply called the {\em cells of the tableau}. The path defining the tableau can be called the {\em shape} of the tableau.
\item a set of blue and red dots in the cells of the tableau such that:
\begin{enumerate}
\item there is no dot below a red dot;
\item there is no dot on the left of a blue dot;
\item any cell of the tableau is either below a red dot, or on the left of a blue dot.
\end{enumerate}
\end{itemize}

Let us give an example of a Catalan alternative tableau of size 23.

\centerline{\epsffile{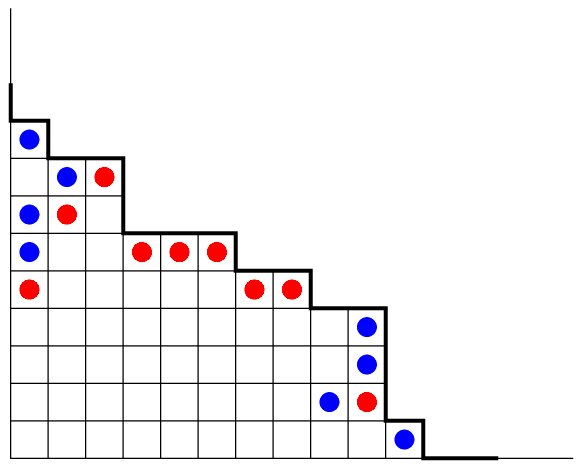}}

It is possible to directly check that Catalan alternative tableaux are enumerated by Catalan numbers (whence their name), but we shall use the following proposition, more adapted to our context.

\begin{prop}\label{bij}
The Catalan alternative tableaux of size $n-1$ are in bijection with binary trees with $n$ internal nodes.
\end{prop}

\proof
We refer to \cite{viennot} (algorithm 2.2) for a formal proof and give here only the idea of the construction. In fact, in that paper, the algorithm was given in term of ``Catalan permutation tableaux'', the subclass of permutation tableaux corresponding to Catalan alternative tableaux, and discussed in \cite{SW}.

We start with Catalan alternative tableau of size $n-1$ and rotate it:

\centerline{\epsffile{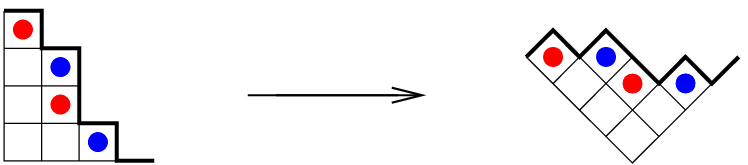}}

The thick line represents the tree under construction. Recursively, for any ``Up-Down'' pattern 

\centerline{\epsffile{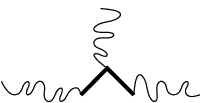}}

in the thick line, we operate a {\em shift}, and two cases are to be distinguished:

\begin{itemize}
\item if the corresponding corner in the tableau is blue: 

\centerline{\epsffile{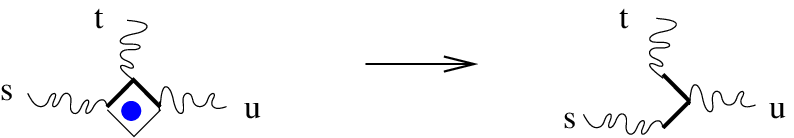}}

and we erase the row of the corner in the tableau;
\item if the corresponding corner in the tableau is red: 

\centerline{\epsffile{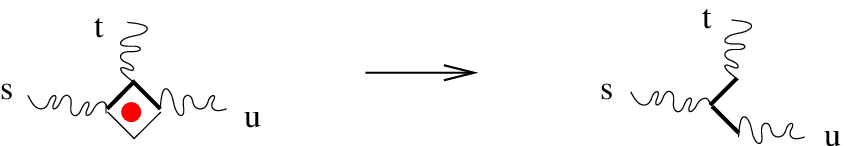}}

and we erase the column of the corner in the tableau.
\end{itemize}

On the example, we obtain:

\centerline{\epsffile{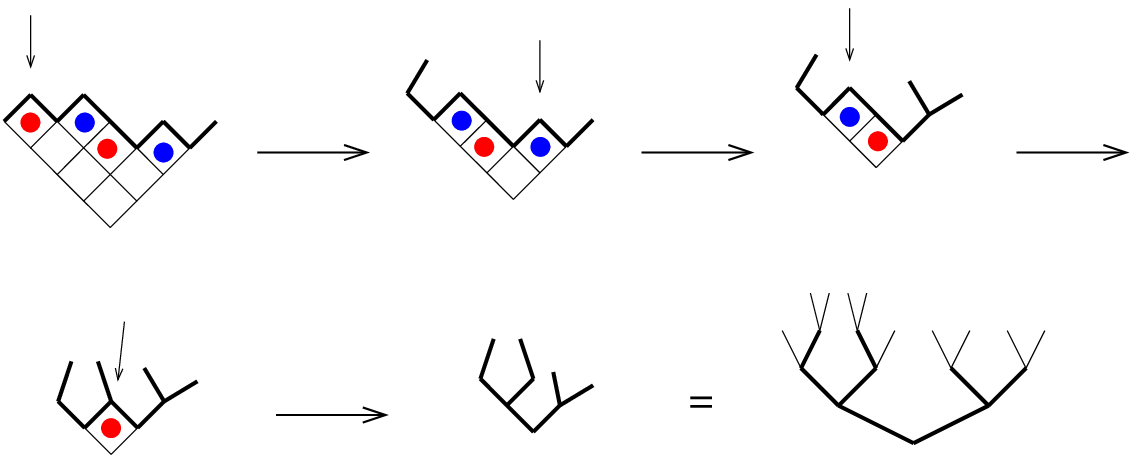}}

It is not difficult to verify that this construction is a bijection.
\endproof

For a permutation $\sigma\in S_n$, its {\em Up-Down sequence} (\cf \cite{viennot2}) is the vector $Q(\sigma)=(q_1,\dots,q_{n-1})\in \{-1,+1\}^{n-1}$ such that
$$q_i=+1\ {\rm iff}\ \sigma_{i+1}>\sigma_i.$$

It is clear that for any tree $T$, all the $\sigma$ in $Z_T$ have the same Up-Down sequence, which we may call the Up-Down sequence of $T$, also called {\em canopy} of the binary tree $T$ in \cite{viennot2}.

Now we may view the shape of a Catalan alternative tableau of size $n-1$ as a vector in $\{-1,+1\}^{n-1}$ (horizontal steps correspond to ``-1'' entries and vertical steps to $+1$ entries).

We have the following property:

\begin{prop}\label{prop1}
The shape of the tableau associated to a tree $T$ through the bijection described in Proposition \ref{bij} is the Up-Down sequence of $T$, as well as the common Up-Down sequence of any permutation $\sigma$ in $Z_T$.
\end{prop}

Now the algorithm described above may be extended to labelled trees: we may put $n$ labels on the shape of a Catalan alternative tableau of size $n-1$:

\centerline{\epsffile{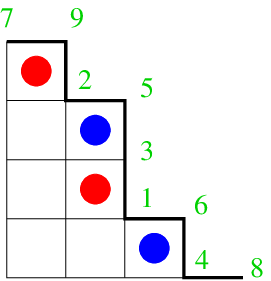}}

If we keep the labels of the nodes when we apply the algorithm to get a tree from the tableau, we obtain a labelled tree

\centerline{\epsffile{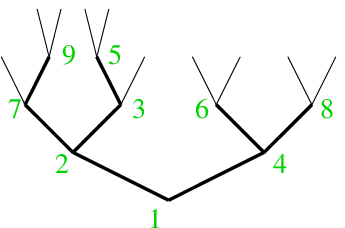}}

In the previous example, we may say that the tableau was labelled with the permutation $792531648$. As a consequence of the bijection, we get the following property.

\begin{prop}\label{prop2}
Let $\sigma$ be a permutation of size $n$. We may label a tableau $C$ with $\sigma$, then apply the bijective algorithm. The labelled tree that we obtain is increasing if and only if:
\begin{itemize}
\item the shape of $C$ is the Up-Down sequence of $\sigma$;
\item the position of the red and blue dots in $C$ is the only one which gives the binary tree $\Psi(\sigma)$.
\end{itemize}
\end{prop}

\section{The product of trees through Catalan alternative tableaux}

Now we come to the main result of this work.
\begin{theo}\label{main}
Let $T_1$ and $T_2$ be two binary trees. Their product in the Loday-Ronco algebra
$$T_1 * T_2=\sum T$$
is given by taking the sum over the trees $T$ associated to Catalan alternative tableaux in the union

\centerline{\epsffile{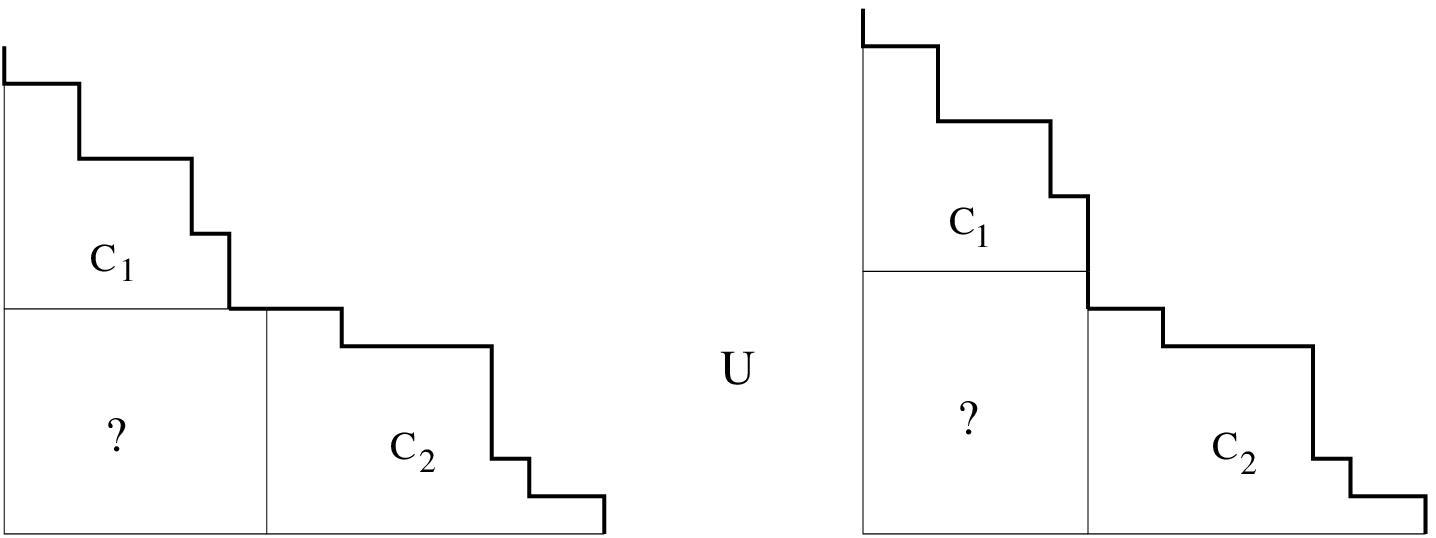}}
where $C_1$ and $C_2$ are the Catalan alternative tableaux associated respectively to $T_1$ and $T_2$, and the question mark ({\rm?}) represent any (valid) placement of (red and blue) dots in the rectangles.
\end{theo}

\proof
By definition of $\Psi^*$, we have:
\begin{equation}\label{eq1}
\Psi^*(T_1 * T_2)=\Psi^*(T_1) * \Psi^*(T_2)=\sum_{\sigma_1\in Z_{T_1}} * \sum_{\sigma_2\in Z_{T_2}}\sigma_2=\sum_{\sigma\in\S} \sigma.
\end{equation}

Let $\sigma$ be an element of $\S$. By definition of the product $*$ in the Malvenuto-Reutenauer algebra, $\sigma$ is of the form: $\sigma=\tau_1.\tau_2$ (concatenation) with the letters appearing in $\tau_1$ and $\tau_2$ form a partition of $\{1,\dots,n\}$, and 
\begin{equation}\label{eq2}
\Psi(\tau_1)=T_1 \ \ {\rm and}\ \ \Psi(\tau_2)=T_2.
\end{equation} 

Thus if $\Psi(\sigma)=T$,the Up-Down sequence of $T$ $Q(T)$ is either $Q(T_1) Up Q(T_2)$ or $Q(T_1) Down Q(T_2)$. Hence the form of the Catalan alternative tableau $C$ associated to $T$ is one of the two given in the Theorem \ref{main}. We label the shape of $C$ with the entries of $\sigma$. The red and blue dots in $C$ have to be placed in a position such that by applying the bijective algorithm, we obtain an increasing binary tree. But if we apply the algorithm to the part of $C$ that carries the entries of $\tau_1$ (respectively $\tau_2$), Propositions \ref{bij} and \ref{prop2} imply that the (red and blue) dots of $C$ in the corresponding subparts of $C$ have to be placed in the same configuration than in $C_1$ (respectively $C_2$). Thie implies that $C$ has the required form.

\medskip 
Conversely, let $T$ be a tableau of the form described in the Theorem \ref{main}, and $\sigma\in Z_T$.
By cutting $\sigma$ in two parts $u_1$ and $u_2$ of lengths the sizes of $C_1$ and $C_2$, we may write: 
$\sigma=u_1.u_2$ with Std$(u_1)=\tau_1$ and Std$(u_2)=\tau_2$.

It is again a simple application of Propositions \ref{bij} and \ref{prop2} that we have: $\Psi(\tau_1)=T_1$ and $\Psi(\tau_2)=T_2$, which was to be proved to complete the proof of Theorem \ref{main}.
\endproof

\section{The $\sha$ product of binary trees}

In light of Theorem \ref{main}, it seems natural to introduce a new product on $k[Y_\infty]$ as follows. 

\begin{defi}
We define the $\sha$ product of two binary trees $T_1$ and $T_2$, associated respectively to Catalan alternative tableaux $C_1$ and $C_2$ by:
$$T_1 \sha T_2=\sum T$$
where the sum is taken over the trees $T$ associated to Catalan alternative tableaux in the set:

\centerline{\epsffile{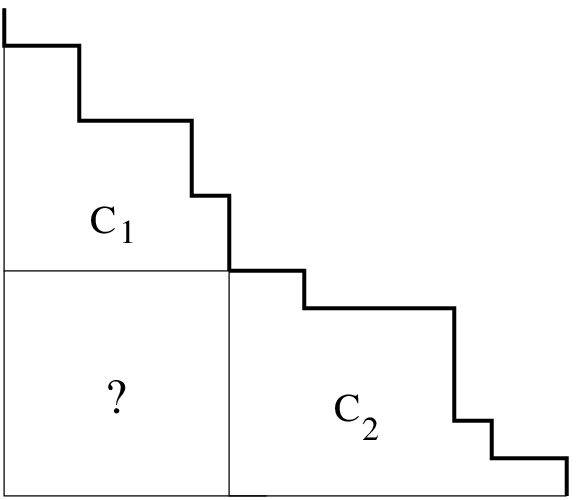}}
and the question mark ({\rm?}) represent any (valid) placement of (red and blue) dots in the rectangles.
\end{defi}

It is clear that this defines an associative product on $k[Y_\infty]$. It is worth to note that for $T_1\in Y_k$ and $T_2\in Y_l$, then the product $T_1 \sha T_2$ is in $Y_{k+l-1}$ (in this case the number of internal edges is preserved).

We give below an example of this product, that should be checked by the reader.

$$\epsffile{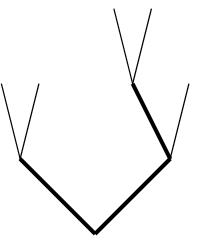}\sha\epsffile{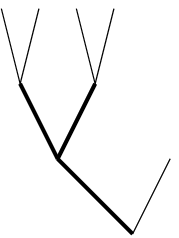}=\epsffile{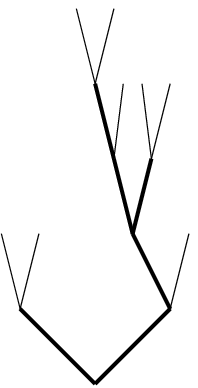}+\epsffile{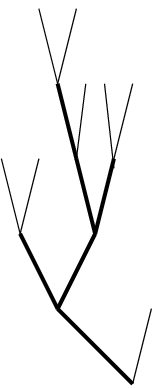}+\epsffile{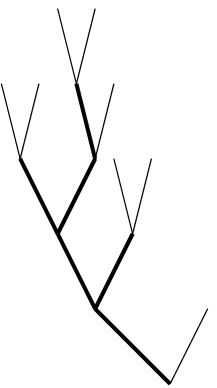}$$

\bigskip \noindent {\bf Acknowledgement.}
The authors sincerely thank F. Chapoton for fruitful discussions.


\end{document}